\documentclass[12pt]{amsart}
\usepackage{amscd,amssymb,xypic}
\input xy
\xyoption{all} 
\newcommand{\bdes}{\begin{description}}
\newcommand{\edes}{\end{description}}
\newcommand{\beqn}{\begin{equation}}
\newcommand{\eeqn}{\end{equation}}

\newcommand{\PP }{{\mathbb P}}
\newcommand{\QQ }{{\mathbb Q}}
\newcommand{\CC }{{\mathbb C}}
\newcommand{\HH }{{\mathbb H}}
\newcommand{\ZZ }{{\mathbb Z}}

\newcommand{\TT }{{\mathbb T}}
\newcommand{\hb }{{\hbar}}
\newcommand{\cT}{{\mathcal T}}
\newcommand{\co}{{\mathcal O}}
\newcommand{\cl}{{\mathcal L}}
\newcommand{\ct}{{c_{\text{top}}}}
\newcommand{\cf}{{\text{ft}}}

\newtheorem{proposition}{Proposition}[subsection]
\newtheorem{lemma}{Lemma}[subsection]

\newtheorem{conjecture}{Conjecture}[subsection]
\newtheorem{def/th}{Definition/Theorem}[subsection]

\begin{document}
\title{A mirror conjecture for projective bundles}
\thanks{1991 Mathematics Subject Classification. Primary 14N35. Secondary 14L30.}
\author{Artur Elezi}
\begin{abstract}
We propose, motivate and give evidence for a relation between the
$\mathcal D$-modules of the quantum cohomology of a smooth complex
projective manifold $X$ and a projective bundle $\PP(\oplus L_i)$
over $X$.
\end{abstract}
\maketitle
\section{\bf Introduction and results}
\noindent Let $Y$ be a projective manifold. Denote by
$Y_{k,\beta}$ the moduli stack of rational stable maps of class
$\beta\in H_2(Y,\ZZ)$ with $k$-markings \cite {[FP]} and
$[Y_{k,\beta}]$ its virtual fundamental class
\cite{[BF]},\cite{[LT]}. Throughout this paper we will be
interested mainly in $k=1$. Recall the following features
\begin{itemize}
\item $e:Y_{1,\beta}\rightarrow Y$ - the evaluation map. \item
$\psi$ - the first chern class of the cotangent line bundle on
$Y_{1,\beta}$. \item $\cf: Y_{1,\beta}\rightarrow Y_{0,\beta}$ -
the forgetful morphism.
\end{itemize}
\noindent Let $\hbar$ be a formal variable and \beqn \label{eq: J}
J_{\beta}(Y):=e_*\left(\frac{[Y_{1,\beta}]}{\hbar(\hbar-\psi)}\right)=\sum_{k=0}^{\infty}
\frac{1}{\hbar^{2+k}}e_*(\psi^k\cap [Y_{1,\beta}]). \eeqn The sum
is finite for dimension reasons. Let $p=\{p_1,p_2,...,p_k\}$ be a
nef basis of $H^2(Y,\QQ)$. For $t=(t_0,t_1,...,t_k)$ let
$$tp:=t_0+\sum_{i=1}^{k}t_ip_i.$$ The $\mathcal D$-module
for the quantum cohomology of $Y$ is generated by \cite{[G2]}
$$J_Y=\exp\left(\frac{tp}{\hb}\right)\sum_{\beta\in
H_2(Y,\ZZ)}q^{\beta}J_{\beta}(Y)$$ where we use the convention
$J_0=1$. The generator $J_Y$ encodes {\it all} of the one marking
Gromov-Witten invariants and gravitational descendants of $Y$.

\noindent For any ring $\mathcal A$, the formal completion of
$\mathcal A$ along the semigroup $MY$ of the rational curves of
$Y$ is defined to be
\begin{eqnarray}
\mathcal A[[q^{\beta}]]:=\{\sum_{\beta\in
\mathrm{MY}}a_{\beta}q^{\beta}, & a_{\beta}\in \mathcal A, &
\beta-\text{effective}\}.
\end{eqnarray}
where $\beta\in H_2(Y,\ZZ)$ is {\it effective} if it is a positive
linear combination of rational curves. This new ring behaves like
a power series since for each $\beta$, the set of $\alpha$ such
that $\alpha$ and $\beta- \alpha$ are both effective is finite. We
may identify $q^{\beta}$ with $q_1^{d_1}\cdot ...\cdot
q_k^{d_k}=\exp(t_1d_1+...+t_kd_k)$ where $\{d_1,d_2,...,d_k\}$ are
the coordinates of $\beta$ relative to the dual basis of
$\{p_1,...,p_k\}$. We view $J_Y$ as an element of
$H^*(Y,\QQ)[t][[q^\beta]]$.

\noindent For toric varieties $J_Y$ is related to an explicit
hypergeometric series $I_Y$ via a change of variables
(\cite{[G1]}, \cite{[LLY3]}). Furthermore, if $Y$ is Fano then the
change of variables is trivial, i.e. $J_Y=I_Y$ thus completely
determining the one point Gromov-Witten invariants and
gravitational descendants of $Y$. This fact was known even earlier
for $Y=\PP^n$ (see for example \cite{[G3]}), which is the case of
relevance for us. An elementary proof of its equivariant version
has been presented in part $3$ of section $3$. Its nonequivariant
limit yields $J_{\PP^n}=I_{\PP^n}$. (An alternative proof in terms
of related torus actions and equivariant considerations follows as
a trivial application from \cite{[Be]}, \cite{[G2]},
\cite{[LLY1]}.) We seek to extend this result in the case of a
projective bundle.

\noindent Let $X$ be a projective manifold. Following
Grothendieck's notation, let $\pi: {\PP}(V)=\PP(\oplus_{j=0}^{n}
L_j)\rightarrow X$ be the projective bundle of hyperplanes of a
vector bundle $V$. Assume that $L_0=\mathcal O_X$. The
$H^*X$-module $H^*{\PP}(V)$ is generated by $z:=c_1(\mathcal
O_{\PP(V)}(1))$ with the relation
$$z \prod_{i=1}^{n}(z-c_1(L_i))=0.$$ Let $s_i: X\rightarrow \PP(V)$ be the section
of $\pi$ determined by the $i$-th summand of $V$ and
$X_i:=s_i(X)$. Then $\mathcal O_{\PP(V)}(1)|_{X_i}\simeq L_i$. Let
$\{p_1,...,p_k\}$ be a nef basis of $H^2(X,\QQ)$.
\begin{lemma} If the line bundles $L_i, i=1,...,n$ are nef then

\noindent (a) $\{p_1,...,p_k,z\}$ is a nef basis of
$H^2(\PP(V),\QQ)$.

\noindent (b) The Mori cones of $X$ and $\PP(V)$ are related via
$$ M\PP(V)=MX\oplus \ZZ_{\geq 0}\cdot [\text{l}]$$
where $[\text{l}]$ is the class of a line in the fiber of $\pi$.
\end{lemma}
\begin{proof} Consider the fibration
\[
\xymatrix
{{\PP^n} \ar[r] & {\PP(V)} \ar[d]_{\pi} \\
& X.}
\]
\noindent There is a first quadrant homology Leray-Serre spectral
sequence with
$$E^2_{p,q}=H_p(X,H_q(\PP^n,\ZZ))$$ which abuts to
$H_{p+q}(\PP(V),\ZZ)$. The diagonal $p+q=2$ and the fact that
$H_1(\PP^n,\ZZ)=0$ yield a short exact sequence $$0\rightarrow
H_2{{\PP}^n}\rightarrow H_2{{\PP}(V)}\stackrel{\pi_*}{\rightarrow}
H_2{X}\rightarrow 0.$$ The pushforward
$${s_0}_*:H_2(X)\rightarrow H_2(\PP(V))$$ splits the above sequence.
Note that for any $\beta \in H_2(X,\QQ)$
\begin{equation}
{s_0}_*(\beta)\cdot z=0. \label{eq: i}
\end{equation}
Let $C$ be a curve in $\PP(V)$ and $f:C\rightarrow \PP(V)$ the
inclusion map. There is a surjection
\begin{equation}
\pi^*(V)\rightarrow {\mathcal O}_{\PP(V)}(1)\rightarrow 0.
\end{equation}
Restricting this sequence to $C$ we obtain
\begin{equation}
\oplus_{i}f^*(L_i)\rightarrow {\mathcal O}_C(z\cdot C)\rightarrow
0.
\end{equation}
Since deg$f^*(L_i)\geq 0$ for all i, we obtain that $z\cdot C \geq
0$. Hence $z$ is nef.  Now if we let $\beta=\pi_*(C)$ then
$\pi_*(C-{s_0}_*(\beta))=0$. Therefore $C-{s_0}_*(\beta)=\nu\cdot
[l]$. But $$\nu\cdot (C\cdot [l])=z\cdot (C-{s_0}_*(\beta))=z\cdot
C.$$ Since $z$ is nef we get $\nu\geq 0$.
\end{proof}
\noindent As seen in the above lemma, if $C\in \PP(V)$ is a curve,
there exists a unique pair $(\nu\geq 0,\beta\in MX$) such that
$[C]=\nu[l]+\beta$. We will identify the homology class $[C]$ with
$(\nu,\beta)$. Note that the generator $J_{\PP(V)}$ is an element
of $H^*(\PP(V),\QQ)[t,t_{k+1}][[q_1^{\nu},q_2^\beta]].$

\noindent To understand the the relation between the generators
$J_{\PP(V)}$ and $J_X$ we first consider the case when the
projective bundle is trivial $\PP(V)=\PP^n\times X$. From now on
we will suppress $X$ and use $J_{\beta}$ instead of
$J_{\beta}(X)$. Denote the homology class of a rational curve in
$\PP^n$ by its degree $\nu$ and the hyperplane class by $H$. Then
$$J_{\PP^n\times X}=J_{\PP^n}\times
J_X=\exp\left(\frac{tp+t_{k+1}H}{\hb}\right)\sum_{\nu \geq
0,\beta} q_1^{\nu}q_2^{\beta}
\frac{1}{\prod_{m=1}^{\nu}(H+m\hb)^{n+1}}J_{\beta}.$$

\noindent For a general vector bundle $V$ we propose twisting the
coefficient in front of $J_{\beta}$ by the chern classes of $V$.
More precisely, let $v_i^{\beta}:=\beta\cdot c_1(L_i)$ and define
the ``twisting'' factor
$$\cT_{\nu,\beta}:=\prod_{i=0}^n{\frac{\prod_{m=-\infty}^{0}(z-c_1(
L_i)+m\hb)}{ \prod_{m=-\infty}^{\nu-v_i^{\beta}}(z-c_1(
L_i)+m{\hb})}}.$$ Let $I_{\nu,\beta}:=\cT_{\nu,\beta}\cdot
{\pi}^*J_{\beta}$ where $\pi^*$ is the flat pull back and define a
``twisted" hypergeometric series for the projective bundle
$\PP(V)$:
$$I_{{\PP}(V)}:=\exp\left(\frac{tp+t_{k+1}z}{\hb}\right)\cdot\sum_{\nu,\beta}{{q_1}^{\nu}}{{q_2}^{\beta}}I_{\nu,\beta}
\in H^*(\PP(V),\QQ)[t,t_{k+1}][[q_1^{\nu},q_2^\beta]].$$ \noindent
In the next section we motivate and propose the following
\begin{conjecture}
Let $L_i, ~i=1,...,n$ be nef line bundles such that $-K_X-c_1(V)$
is ample. Then $J_{\PP(V)}=I_{\PP(V)}.$
\end{conjecture}
\noindent There is much evidence for this conjecture. First the
extremal cases: on one end, the pure fiber case, i.e. the equality
of the $\beta=0$ terms is known to be true; on the other end we
prove the equality of the $d=0$ terms when all the line bundles
$L_i$ are ample. We also show that this conjecture is consistent
with the Quantum Lefshetz Principle applied to all the complete
intersections $\displaystyle{\prod_{j\neq i}(z-c_1(L_i))}$.
Finally, when $X$ is a toric variety and $L_i$ are toric line
bundles the projective bundle $\PP(V)$ is also a toric variety and
the validity of the conjecture follows from the toric mirror
theorem (\cite{[G1]},\cite{[LLY3]}).

\section{\bf Motivation}
\noindent In this section we motivate the conjecture by comparing
$J_{\PP(V)}$ and $I_{\PP(V)}$. Recall first that there is a
natural grading in the quantum cohomology of $\PP(V)$
(\cite{[G1]},\cite{[G2]})
\begin{itemize}
\item $\text {deg}~q_1:=n+1$ \item
$\text{deg}~q_2^{\beta}:=\beta\cdot (-K_X-c_1(V))$ \item
$\text{deg}(\alpha)=\text{codim}_{\CC}~(\alpha)$ for any pure
cycle $\alpha\in H_*(\PP(V),\ZZ)$.
\end{itemize}
\begin{lemma} \label{lemma: homogeneity}
Both $J_{\PP(V)}$ and $I_{\PP(V)}$ are homogeneous of degree zero.
\end{lemma}
\begin{proof}
From the exact sequence
\[0\rightarrow {\mathcal
O}_{\PP(V)}\rightarrow \pi^*V^*(1)\rightarrow
T_{\PP(V)}\rightarrow \pi^*T_X\rightarrow 0\] we find that
\begin{equation}
K_{\PP(V)}=\pi^*K_{X}+\pi^*c_1(V)-(n+1)z.
\end{equation}
Multiplying by a curve class $(\nu,\beta)$ in $H_2\PP(V)$ we
obtain
\[-(\nu,\beta)\cdot K_{\PP(V)}=(n+1)\nu-\beta\cdot K_X-\beta\cdot c_1(V)=\text{deg}~(q_1^\nu)\cdot\text{deg}~(q^{\beta}).\]
Hence
\[
\text{dim}~e_*([\PP(V)_{1,(\nu,\beta)}])=\text{deg}~(q_1^\nu)\cdot
\text{deg}~(q^{\beta})+\text{dim}~\PP(V)-2
\]
Therefore
\[
\text{deg}~e_*(\psi^k\cap
[\PP(V)_{1,(\nu,\beta)}])=2+k-\text{deg}~(q_1^\nu)\cdot\text{deg}~(q^{\beta}).
\]
It follows from the presentation (\ref{eq: J}) that the degree of
$J_{\nu,\beta}$ is
$$-\text{deg}~(q_1^{\nu})\cdot\text{deg}~(q^{\beta}),$$ hence
$\sum q_1^{\nu}q^{\beta}J_{\nu,\beta}(\PP(V))$ is homogeneous of
degree zero.

\noindent We turn our attention to the generator $I_{\PP(V)}$.
Notice that
$$\text{deg}~e_*(\psi^k\cap [X_{1,\beta}])=\text{dim}~X-(\text{dim}~[X_{1,\beta} ]-k)=\beta\cdot K_X+2+k$$ and so
\[\text{deg}~(\pi^*(J_{\beta}))=\text{deg}~J_{\beta}=\beta\cdot K_X.\]
Now
\[
\text{deg}~{\bar I}_{\nu,\beta}=-\beta\cdot c_1(V)-(n+1)\nu.
\]
It follows that
$$\text{deg}~I_{\nu,\beta}=\beta\cdot c_1(V)-(n+1)\nu+\beta\cdot K_X=-\text{deg}~(q_1^\nu)\cdot\text{deg}~(q^{\beta})$$
and the lemma is proved.
\end{proof}
For a trivial bundle $V$, we showed that $J_{\PP(V)}=I_{\PP(V)}$
using the corresponding identity for projective spaces. From the
general framework of mirror theorems, it is reasonable to expect
this identity to hold for general $V$, if the $\hb$-asymptotic
expansions of $J_{\PP(V)}$ and $I_{\PP(V)}$ coincide in orders
$\hb^0$ and $\hb^{-1}$ (\cite{[G2]}). By the definition
$$\sum_{\nu,\beta}{{q_1}^{\nu}}{{q_2}^{\beta}}
J_{\nu,\beta}=1+o({\hbar}^{-1}).$$ We now show that, under the
assumptions of the conjecture, $I_{\PP(V)}$ has the same
asymptotic expansion. For $0\neq (\nu,\beta)\in M\PP(V)$, let
$l_{\nu,\beta}$ be the cardinality of the set
$$S_{\nu,\beta}=\{i:\nu-v_i^{\beta}<0\}.$$
A simple algebra shows that
$$\displaystyle{\bar
I_{\nu,\beta}=\frac{a_{\nu,\beta}}{\hbar^{l_{\nu,\beta}+(n+1)\nu-c_1(V)\cdot
{\beta} }}\cdot P(\hbar^{-1}}).$$ Here $a_{\nu,\beta}\in
H_*(\PP(V))$ is homogeneous of degree $l_{\nu,\beta}$ and
$P(\hbar^{-1})$ is a power series of $\hb^{-1}$. On the other
hand, it follows from $(1)$ that
$$\pi^*(J_{\beta})=\frac{a_{\beta}}{\hbar^{n_{\beta}}}\cdot Q(\hbar^{-1})$$
where $n_{\beta}=\mathrm{max}(2,-\beta\cdot K_X)$ and
$Q(\hbar^{-1})$ is also a power series of $\hbar^{-1}$. It follows
that $I_{\nu,\beta}={\bar I}_{\nu,\beta}\cdot \pi^*(J_{\beta})$ is
a power series of $\hbar^{-1}$ with leading term
$$\frac{A_{\nu,\beta}}{\hbar^{l_{\nu,\beta}+(n+1)\nu-\beta\cdot c_1(V)+n_{\beta}}}$$ for
every $0\neq (\nu,\beta)\in M\PP(V)$. Let
$n_{\nu,\beta}:=l_{\nu,\beta}+(n+1)\nu-\beta\cdot
c_1(V)+n_{\beta}.$ If $\nu> 0$ then $n_{\nu,\beta}>(n+1)\nu>1.$ If
$\beta \cdot c_1(V)=0$ then $v_i^{\beta}=0, i=0,1,...,n$ hence
$l_{0,\beta}=0$ and $n_{0,\beta}=n_{\beta}\geq 2.$ Finally if
$\beta \cdot c_1(V)>0$ then $l_{0,\beta}>0$ hence $n_{0,\beta}>
(-K_X-c_1(V))\cdot \beta\geq 1$. It follows that
$I_{\nu,\beta}=o({\hbar}^{-1})$ for any $0\neq (\nu,\beta)\in
M\PP(V).$
\section{\bf Evidence for the conjecture}
\noindent {\bf 1. The conjecture holds for toric varieties.}
Indeed, assume $X$ is a toric manifold and let $\Sigma\subset
\ZZ^m$ be its fan. Denote by $b_1,...,b_r$ its one dimensional
cones. If $L_i, i=0,1,...,n$ are toric line bundles then the
projective bundle $\PP(\oplus_{i=0}^n L_i)\rightarrow X$ is also a
toric variety and there is a canonical way to obtain its fan
\cite{[O]}. Let $\ZZ^n$ be a new lattice with basis
$\{f_1,...,f_n\}$. The edges $b_1,...,b_r$ of $\Sigma$ are lifted
to new edges $B_1,B_2,...,B_r$ in $\ZZ^m\oplus \ZZ^n$ and
subsequently $\Sigma$ is lifted in a new fan $\Sigma_1$ in the
obvious way. Let $\Sigma_2\subset 0\oplus \ZZ^n$ be the fan of
$\PP^n$ with edges $F_0=-\sum_{i=1}^n f_i,F_1=f_1,...,F_n=f_n$.
The canonical fan associated to $\PP(V)$ consists of the cones
$\sigma_1+ \sigma_2$ where $\sigma_1, \sigma_2$ are cones in
$\Sigma_1,\Sigma_2$. Let $N=r+n+1$. The edges $F_i$ in the
canonical fan of the projective bundle, correspond to the divisors
$z-c_1(L_i)$ where $z:=c_1(\co_{\PP(V)}(1))$ while $B_i$
correspond to the pullback of the base divisors \cite{[M]}. Since
$L_i, i=0,1,...,n$ are nef and $-K_X-\sum_i L_i$ is ample it
follows that $-K_{\PP(V)}$ and $-K_X$ are also ample. Hence from
\cite{[G1]}
$$J_{\PP(V)}=\exp\left(\frac{tp+t_{k+1}z}{\hb}\right)
\cdot\sum_{(\nu,\beta)\in
MP(V)}{{q_1}^{\nu}}{{q_2}^{\beta}}\tau_{\nu,\beta}\prod_i
\frac{\prod_{m=-\infty}^{0}(b_i+m\hb)}{\prod_{m=-\infty}^{b_i\cdot
\beta}(b_i+m\hb)}$$ and
$$J_X=\exp\left(\frac{tp}{\hb}\right)\cdot\sum_{\beta\in MX}{{q_2}^{\beta}}\prod_i
\frac{\prod_{m=-\infty}^{0}(b_i+m\hb)}{\prod_{m=-\infty}^{b_i\cdot
\beta}(b_i+m\hb)}.$$ The equality $J_{\PP(V)}=I_{\PP(V)}$ follows
readily.

\vspace{.2in} \noindent {\bf 2. Consistency with the Quantum
Hyperplane Section Principle.} Let $W=\oplus_{i=1}^{r}
\cl_i\rightarrow Y$ be a vector bundle over a projective manifold
$Y$. Let $Z=Z(s)$ be the zero locus of a generic section $s$ of
$W$ and $\mu$ the inclusion map. For any $0\neq \beta\in
H_2(Y,\ZZ)$, denote by $W_{\beta}$ the bundle on $\overline
M_{0,1}(Y,\beta)$ whose fiber over a stable map $(C,x_1,f)$ is the
vector space of the sections of $H^0(C,f^*(W))$ that vanish at
$x_1$. The following $H^*(Y)[[q^\beta]]$-valued generator \beqn
J_W=\exp\left(\frac{tp}{\hbar}\right)c_r(W)
\left(1+\sum_{\beta\neq
0}q^{\beta}e_*\left(\frac{\ct(W_{\beta})}{\hbar(\hbar-\psi)}\right)\right).
\eeqn is intrinsically related to the quantum $\mathcal D$-module
of $Z$. Assume that $\mu^*: H^2(Y,\QQ)\rightarrow H^2(Z,\QQ)$ is
surjective. The Gysin map in cohomology extends to a map of
completions $$\mu_*: H^*(Z,\QQ)[t_Z][[q^{\beta}]]\rightarrow
H^*(Y,\QQ)[t_Y][[q^{\alpha}]]$$ via
$\mu^*(q^{\beta})=q^{\mu_*(\beta)}$. It is shown in \cite{[AE]}
that $$\mu_*(J_Z)=J_W.$$ Consider now
\begin{equation}
I_W=\exp\left(\frac{tp}{\hbar}\right)c_r(W)
\left(1+\sum_{\beta\neq 0}q^{\beta}L^W_{\beta}
e_*\left(\frac{[Y_{1\beta}]}{\hbar(\hbar-\psi)} \right)\right)
\end{equation}
where the modifying Lefshetz factor $L^W_{\beta}$ is defined to be
$$\displaystyle{L^W_{\beta}:=\prod_{i=1}^{r}\frac{\prod_{m=-\infty}^{b_i}(c_1(
\cl_i)+m\hbar)}{\prod_{m=-\infty}^{0}(c_1( \cl_i)+m\hbar)}}$$ with
$b_i:=\beta\cdot c_1(\cl_i)$. The Quantum Hyperplane Section
Principle asserts that, under suitable conditions, $J_W$ equals
$I_W$ up to a change of variables of the form
\[t_0'=t_0+f_0(q)\hbar+f(q)\]
\[t'_i=t_i+f_i(q),~i=1,2,...,n\]
where $f_i$'s are homogeneous power series of degree $0$ and $f$
is homogeneous power series of degree $1$ (\cite{[K1]},
\cite{[K2]},\cite{[YPL]}, and for the most general form
\cite{[GiCo]}). Each section $s_i:X\rightarrow X_i\subset \PP(V)$
is the zero locus of a section of the line bundle $\oplus_{j\neq
i}\mathcal O_{\PP(V)}(1)\otimes L_j^*$ in ${\PP}(V)$.
\begin{proposition}
Quantum Hyperplane Section Principle applied appropriately to the
generator $I_{\PP(V)}$ produces the (Gysin image of the) generator
$J_{X_i}$ of the quantum $\mathcal D$-module of $X_i$.
\end{proposition}
\begin{proof} For notational simplicity we focus in the case $i=0$ and let
$\mu_0:X_0\rightarrow \PP(V)$ be the inclusion map. Consider the
following cohomology-valued function
\[I= \exp\left(\frac{tp+t_{k+1}z}{\hb}\right) \prod_{i=1}^{n}(z-c_1(L_i))
 \sum_{d,\beta} {q_1}^d{q_2}^{\beta}\cdot\]
\[
\prod_{k=1}^{n}\frac{\prod_{m=-\infty}^{d-v_k}(z-c_1(
L_k)+m{\hb})}{\prod_{m=-\infty}^{0}(z-c_1(
L_k)+m{\hb})}\prod_{i=0}^{n}\frac{\prod_{m=-\infty}^{0}(z-c_1(
L_i)+m\hb)}{\prod_{m=-\infty}^{d-v_i}(z-c_1(
L_i)+m{\hb})}{\pi}^*J_{\beta}=
\]
\[
\exp\left(\frac{tp+t_{k+1}z}{\hb}\right)
\prod_{i=1}^{n}(z-c_1(L_i))
 \sum_{d,\beta} {q_1}^d{q_2}^{\beta}\cdot\frac{1}{\prod_{m=1}^{d}(z+m{\hb})}{\pi}^*J_{\beta}
\]
Expanding
$$\frac{1}{\prod_{m=1}^{d}(z+m{\hb})}$$
as a power series in $z$ and using the relation
$$z\prod_{i=1}^{n}(z-c_1(L_i))=0,$$  we find that
\[I=\prod_{i=1}^{n}(z-c_1(L_i))\text{exp}\left(\frac{tp}{\hbar}\right)
\exp\left(\frac{q_1}{\hb}\right)\sum_{\beta}q_2^{\beta}\pi^*(J_{\beta}).\]
Let ${t'}_0=t_0+q_1$. With this new variables we have
${\mu_0}_*(J_{X_0})=I$. The proposition is proven.
\end{proof}
\vspace{.2in} \noindent {\bf 3. Extremal Case I: Pure fiber}. In
the case where $\beta=0$ our conjecture reads
$$e_*\left(\frac{[\PP(V)_{1,(\nu,0)}]}{\hb(\hb-\psi)}\right)=\prod_{i=0}^n{\frac{1}{\prod_{m=1}^{\nu}(z-c_1(
L_i)+m{\hb})}}.$$ This is true. In fact the restriction of
$J_{\PP(V)}$ to $\beta=0$ is obtained from the relative
Gromov-Witten theory of $\PP(V)$ over $X$, which can be defined
via intersection theory in the spaces of holomorphic curves to the
fibers of $\pi$. If we let $\TT=(\CC^*)^{n+1}$ act diagonally on
$\PP^n$ then the relative GW theory of the $\PP^n$ bundle over
$B\TT$ associated with the universal bundle $E\TT\mapsto B\TT$ is
precisely the $\TT$-equivariant GW theory of $\PP^n$
(\cite{[GK]}). Furthermore, if $f:X \mapsto B\TT$ induces the
principal bundle associated with $\PP(V)$, then the
$\TT$-equivariant theory of $\PP^n$ pulls back via $f$ to the
relative GW theory of $\PP(V)$ over $X$. It follows that the
restriction of $J_{\PP(V)}$ to $\beta=0$ is obtained by
substituting $c_1(L_i)$ for the generators $\lambda_i$ of
$H^*(B\TT)=\QQ[\lambda_0,...,\lambda_n]$ in the $\TT$-equivariant
generator $J^{\TT}_{\PP^n}$. Now recall that $J^{\TT}_{\PP^n}$ and
its derivatives form a fundamental solution matrix to the quantum
differential equation
$$\hb q\frac{dS}{dq}=H\ast S,$$ where we use the same letter $H$ to denote the
equivariant hyperplane class and $\ast$ is the product in the
small equivariant quantum cohomology algebra of $\PP^n$ (this is
essentially due to Dijkgraaf-Witten and Dubrovin \cite{[D]}). For
dimension reasons, if $k<n$ we have $H\ast H^k=H^{k+1}$, where the
powers of $H$ are the classical ones. The product $H\ast H^n$ is
found from the only possible graded deformation
$$\prod_{i=0}^n (H-\lambda_i)=xq$$ of the classical relation (where
$x=0$). The coefficient $x$ is equal to $1$ (two points in $\PP^n$
determine a unique line). The small equivariant quantum cohomology
of $\PP^n$ is
 $$\QQ[H,\lambda_0,...,\lambda_n,q]/(\prod_{i=0}^n(H-\lambda_i)-q)$$
(see also \cite{[AS]}, \cite{[K3]}.) The relations imply that
$J^{\TT}_{\PP^n}$ is annihilated by the operator
$$D=\prod_{i=0}^n (\hb q\frac{d}{dq}-\lambda_i)-q.$$ One can see that $J^{\TT}_{\PP^n}$
is the only solution of the above equation that has the form
$$\exp\left(\frac{t_0+H\ln q}{\hb}\right)(1+q\sum_{\nu\geq
0}a_{\nu}q^{\nu}),$$ where $a_{\nu}\in
H^*_{\TT}(\PP^n)[\hb^{-1}]$. But
$$I^{\TT}_{\PP^n}=\exp\left(\frac{t_0+H\ln q}{\hb}\right)\sum_{\nu
\geq 0} q^{\nu} \prod_{i=0}^n
\frac{1}{\prod_{m=1}^{\nu}(H-\lambda_i+m\hb)}$$ has also the above
form and a direct calculation shows that it is annihilated by the
operator $D$. Hence
\begin{equation}
\label{equiv} J^{\TT}_{\PP^n}=I^{\TT}_{\PP^n}
\end{equation}
and the pure fiber case follows. (For another approach see
\cite{[Be]}.) We note that the nonequivariant limit (i.e.
$\lambda_i\to 0,~i=0,1,...,n$) of (\ref{equiv}) yields
$J_{\PP^n}=I_{\PP^n}$.

\vspace{.2in}

\noindent {\bf 4. Extremal Case II: No fiber.} This is the case
where $\nu=0$.
\begin{proposition} If the line bundles $L_i, i=1,2,...,n$ are ample then the conjecture holds for $d=0$ i.e.
$$J_{0,\beta}(\PP(V))=\prod_{i=1}^{n}\prod_{m=0}^{-\beta\cdot
c_1(L_i)-1}(z-c_1(L_i)-m\hbar)\pi^*(J_{\beta}).$$
\end{proposition}
\begin{proof} The trivial action of the torus $\TT=(\CC^*)^{n+1}$ on $X$ lifts to a diagonal action on
the fibers of $\PP(V)$. Let $\co(1)_i$ be the pull back of the
$i$-th $\co_{\PP^{\infty}}(1)$ to $X_{\TT}=X\times
\PP^{\infty}\times...\times \PP^{\infty}$ and $\cl_i:=L_i\otimes
\co(1)_i$. Then $\PP(V)_{\TT}=\PP(\oplus_{i=0}^{n}(\cl_i))$ and
$${H_{\TT}}^*(\PP(V)):=H^*(\PP(V)_{\TT})=H^*X[z_{\TT},u_0,...,u_n]/\prod_{i=0}^{n}(z_{\TT}-c_1(L_i)-u_i),$$
where $z_{\TT}:=c_1(\co_{\PP(V)_{\TT}}(1))$ and
$u_i:=c_1(\co(1)_i)$. The connected components of the fixed point
loci are the $n+1$ sections $s_i:X\simeq X_i\subset \PP(V)$ that
correspond to the quotients $V\rightarrow L_i$. Note that
$$s_i^*(z_{\TT})=c_1(L_i)+u_i.$$ The proposition will follow as a nonequivariant
limit (i.e. $u_i=0$) of  the equivariant identity
\[
e_*\left(\frac{[\PP(V)_{1,(0,\beta)}]}{\hbar(\hbar-\psi)}\right)=\prod_{i=1}^{n}\prod_{m=0}^{-\beta\cdot
c_1(L_i)-1}(z^{\TT}-c_1(L_i)-u_i-m\hbar)\pi^*(J_{\beta})
\]
where the left hand side is the equivariant pushforward.

\noindent By the localization theorem (\cite{[AB]}), it suffices
to show that the restrictions to the fixed point components $X_i$
agree. For $i=1,2,...,n$ the restriction of the right hand side
obviously vanishes. The $\TT$-action on $\PP(V)$ induces an action
on $\PP(V)_{1,(0,\beta)}$. If $i\neq 0$, none of the connected
components of $\PP(V)^{\TT}_{1,(0,\beta)}$ is mapped to $X_i$ by
the evaluation map $e$. Hence the restriction of the left side to
$X_i$ vanishes for $i=1,2...,n$.

\noindent We notice that $N_{X_0/\PP(V)}=\oplus_{i=1}^n L_0\otimes
L_i^*\simeq \oplus_{i=1}^{n}L_i^*$. The induced $\TT$ action is
such that
$$e^{\TT}(N_{X_0/\PP(V)})=\prod_{i=1}^{n}(u_0-u_i-c_1(L_i))=\prod_{i=1}^{n}(-c_1^{\TT}(L_i)).$$ The
restriction of the right hand side to $X_0$ is
$$e^{\TT}(N_{X_0/\PP(V)})\prod_{i=1}^{n}\prod_{m=1}^{-\beta\cdot
c_1(L_i)-1}(-c_1^{\TT}(L_i)-m\hbar)J_{\beta}.$$

\noindent Only one fixed point component of
$\PP(V)^{\TT}_{1,(0,\beta)}$ is mapped to $X_0$. It is isomorphic
to $X_{1,\beta}$. Since $f^*(L_i)$ is positive for every
nonconstant map $f:\PP^1\rightarrow X$,
$$\HH^1:=R^1\cf_*e^*(N_{X_0/\PP(V)})=R^1\cf_*e^*(\oplus_{i=1}^{n}L_i^*)$$
is a vector bundle on this fixed point component. A simple
analysis of the obstruction theory shows that the equivariant
Euler class of the normal bundle to this component is inverse to
the equivariant Euler class of $\HH^1$. We apply the localization
theorem for the evaluation map $e:\PP(V)_{1,(0,\beta)}\rightarrow
\PP(V)$ and the fixed point component $X_0\subset \PP(V)$
(\cite{[Be]}, \cite{[LLY2]}) to obtain
$$e_*\left(\frac{\ct(\HH_1)\cap [X_{1,\beta}]}{\hb(\hb-\psi)}\right)=
\frac{1}{e^{\TT}(N_{X_0/\PP(V)})}\cdot
s_0^*\left(e_*\left(\frac{[\PP(V)_{1,(0,\beta)}]}{\hbar(\hbar-\psi)}\right)\right).$$
\noindent Suffices to show that \beqn \label{eq: basic}
e_*\left(\frac{\ct(\HH^1)\cap
[X_{1,\beta}]}{\hbar(\hbar-\psi)}\right)=\prod_{i=1}^{n}\prod_{m=1}^{-\beta\cdot
c_1(L_i)-1}(-c^{\TT}_1(L_i)-m\hbar)J_{\beta}. \eeqn This is an
equivariant form of the concave Quantum Lefshetz Principle for the
diagonal fiber action. The nonequivariant proofs of \cite{[GiCo]}
or \cite{[YPL]} can be trivially modified for this equivariant
version.
\end{proof}

\noindent{\bf Acknowledgement}. The author thanks A. Givental, T.
Graber, S. Katz, B. Kim, Y.-P. Lee, J. Li, K. Liu for very useful
conversations. We are also thankful to the referee for many
valuable suggestions that greatly improved the quality of this
manuscript.

\vspace{+6 pt} \noindent Department of Mathematics and Statistics,

\noindent American University

\noindent 4400 Massachusetts Ave, Washington DC, 20016

\noindent aelezi@american.edu

\end{document}